\documentclass[11pt]{article}
\usepackage{amsfonts,amssymb,amsmath,mathrsfs}
\usepackage{color}
\usepackage[colorlinks, linkcolor=blue, citecolor=blue]{hyperref}
\usepackage{graphicx}
\usepackage[T1]{fontenc}
\usepackage{hyperref}
\hypersetup{
    colorlinks=true,
    linkcolor=blue,
    filecolor=blue,
    urlcolor=blue,
    citecolor=blue
}

\allowdisplaybreaks[4]

\allowdisplaybreaks

\newtheorem{con0}{Theorem}[section]
\newtheorem{thm0}{Theorem}[section]
\newtheorem{exa0}{Theorem}[section]

\newtheorem{con1}[con0]{Condition}
\newtheorem{def1}[thm0]{Definition}
\newtheorem{lem1}[thm0]{Lemma}
\newtheorem{thm1}[thm0]{Theorem}
\newtheorem{cor1}[thm0]{Corollary}
\newtheorem{pro1}[thm0]{Proposition}
\newtheorem{rem1}[thm0]{Remark}
\newtheorem{ass1}[thm0]{Assumption}
\newtheorem{exa1}[exa0]{Example}

\def\blemma{\begin{lem1}}\def\elemma{\end{lem1}}
\def\btheorem{\begin{thm1}}\def\etheorem{\end{thm1}}
\def\bcorollary{\begin{cor1}}\def\ecorollary{\end{cor1}}
\def\bproposition{\begin{pro1}}\def\eproposition{\end{pro1}}
\def\bremark{\begin{rem1}\rm{}\def\eremark{\end{rem1}}}

\def\benumerate{\begin{enumerate}}\def\eenumerate{\end{enumerate}}
\def\bitemize{\begin{itemize}}\def\eitemize{\end{itemize}}\def\itm{\item}

\def\beqlb{\begin{eqnarray}}\def\eeqlb{\end{eqnarray}}
\def\beqnn{\begin{eqnarray*}}\def\eeqnn{\end{eqnarray*}}
\def\beqa{\begin{equation}}\def\eeqa{\end{equation}}
\def\beqaa{\begin{equation*}}\def\eeqaa{\end{equation*}}

\def\eqref#1{{\rm(\ref{#1})}}

\def\ulim{\,\uparrow\!\!\lim}\def\dlim{\,\downarrow\!\!\lim}
\def\ar{\!\!\!&}\def\nnm{\nonumber}\def\ccr{\nnm\\}
\def\qqquad{\qquad\qquad}

\def\bproof{\begin{proof}~}\def\eproof{\qed\end{proof}}
\def\qed{\hfill$\square$\smallskip}

\def\mrm{\mathrm}\def\mbf{\mathbf}
\def\mbb{\mathbb}

\def\d{\mrm{d}}\def\e{\mrm{e}}

\def\blue{\color{blue}}

\begin{document}

\noindent{Published in: Journal of Statistical Physics (2025) 192:98}

\noindent{https://doi.org/10.1007/s10955-025-03480-3}

\bigskip\bigskip

\noindent{\bf\LARGE Asymptotic behavior of the generalized}

\medskip

\noindent{\bf\LARGE Derrida--Retaux recursive model\,\footnote{\,Some errors are corrected in blue colour.}}

\bigskip

\noindent{Zenghu Li and Run Zhang\,\footnote{\,Corresponding author.}}

\medskip{\small

\noindent{\it Laboratory of Mathematics and Complex Systems, School of Mathematical Sciences,}

\noindent{\it Beijing Normal University, Beijing 100875, China}

\noindent{Emails: \tt lizh@bnu.edu.cn, zhangrun@mail.bnu.edu.cn}

\bigskip

\noindent\textbf{Abstract:~} We study the max-type recursive model introduced by Hu and Shi (J. Stat. Phys., 2018), which generalizes the model of Derrida and Retaux (J. Stat. Phys., 2014). The class of geometric-type marginal distributions {\blue is} preserved by the model with {\blue a} geometric offspring distribution. We give some long-time asymptotic expansions of the parameters of the marginal distribution. From the expansions, we derive the asymptotics of the sustainability probability, marginal distribution, first moment and probability generating function.

\bigskip

\noindent\textbf{Keywords:~} Generalized Derrida--Retaux model, geometric offspring distribution, geometric-type distribution, sustainability probability, first moment, generating function, asymptotic behavior

\smallskip

\noindent\textbf{Mathematics Subject Classification:~} 60J05, 82B27

}

\bigskip


\section{Introduction}

\setcounter{equation}{0}

A discrete-time max-type recursive model was introduced by Derrida and Retaux~\cite{DeR14} in their study of the depinning transition in the limit of strong disorder. Let $X_0$ be a random variable taking values in $\mbb{Z}_+:= \{0,1,2,\cdots\}$. By a \textit{Derrida--Retaux process}, or \textit{DR process} for short, we mean a discrete-time stochastic process $(X_n: n\ge 0)$ defined recursively in the distribution sense by
 \beqlb\label{DR0proc}
X_{n+1}\overset{\rm d}= (X_n+\tilde{X}_n-1)_+, \quad n\geq 0,
 \eeqlb
where $(\,z\,)_+:= \max(0,z)$ and $\tilde{X}_n$ is an independent copy of $X_n$. From \eqref{DR0proc} it is easy to see that
 \beqnn
2\mbf{E}(X_n)-1\le \mbf{E}(X_{n+1})\le 2\mbf{E}(X_n).
 \eeqnn
For this model, the \textit{free energy} is defined by the limit:
 \beqnn
F_\infty:= \dlim_{n\to \infty} \,2^{-n}\mbf{E}(X_n)
 =
\ulim_{n\to \infty} \,2^{-n}[\mbf{E}(X_n)-1].
 \eeqnn
The DR process is referred to as \textit{pinned} or \textit{supercritical} if $F_\infty> 0$, and as \textit{unpinned} or \textit{subcritical} if $F_\infty= 0$. The two regimes are separated by the \textit{critical manifold} in a suitable parametric space. There are many interesting problems related to the behavior at or near the criticality.

Let $\mu_n$ denote the distribution of $X_n$ for $n\ge 0$. To discuss the phase transition from the pinned to the unpinned regime, it is convenient to specify the mass of $\mu_0$ at the origin. Consider the decomposition:
 \beqnn
\mu_0= p\delta_0 + (1-p)\vartheta,
 \eeqnn
where $0\le p\le 1$ is a constant and $\vartheta$ is a fixed probability measure carried by $\{1,2,\cdots\}$. Let $F_\infty(p)$ denote the associated free energy. Then $p\mapsto F_\infty(p)$ is a decreasing function on $[0,1]$. Write $p_c\in [0,1]$ for the \textit{critical parameter}, that is,
 \beqnn
p_c= \sup\{p\in [0,1]: F_\infty(p)> 0\}
 \eeqnn
with the convention $\sup\emptyset= 0$. Derrida and Retaux~\cite{DeR14} conjectured that, under the assumption $p_c>0$, there exists some constant $C>0$ such that
 \beqlb\label{DR0conj}
F_\infty(p)= \exp\Big\{\!-\frac{C+o(1)}{(p_c-p)^{1/2}}\Big\}, \quad p\uparrow p_c.
 \eeqlb

A number of variations of the DR process have also been studied. Let $\eta$ be an integrable random variable taking values in $\{1,2,\cdots\}$ with $m:= \mbf{E}(\eta)> 1$. Instead of \eqref{DR0proc}, we can also define a discrete-time max-type process $(Y_n: n\ge 0)$ by the recursive formula
 \beqlb\label{gDR0proc}
Y_{n+1}\overset{\rm d}= (Y_{n,1}+Y_{n,2}+\cdots+Y_{n,\eta}-1)_+,
 \eeqlb
where $\{Y_{n,1}, Y_{n,2}, \ldots\}$ is a sequence of independent copies of $Y_n$ independent of $\eta$. For this model, the \textit{free energy} is defined as the limit
 \beqlb\label{F_infty:=def}
F_\infty:= \lim_{n\to \infty}m^{-n}\mbf{E}(Y_n).
 \eeqlb
The model defined by \eqref{gDR0proc} was studied in Collet et al.~\cite{CEGM84b} as a spin glass model, in Curien and H\'enard \cite{CuH22} and Goldschmidt and Przykucki \cite{GoP19} as a parking scheme, and in Li and Rogers \cite{LiR99} as an iteration function of random variables. It also belongs to the class of max-type recursive models discussed in the survey paper of Aldous and Bandyopadhyay \cite{AiB05}. By the result of Collet et al.~\cite{CEGM84b}, the model exhibits the phase transition given as follows.

\medskip\noindent\textbf{Theorem~A\,} {(Collet et al.~\cite{CEGM84b})\,} \textit{Suppose that $\eta\overset{\rm a.s.}= m\ge 2$ and $\mbf{E}[Y_0m^{Y_0}]< \infty$. Then:}
 \benumerate

\itm[{\rm(i)}] \textit{$\mbf{E}(m^{Y_0})< (m-1)\mbf{E}(Y_0m^{Y_0})$ implies $F_\infty> 0$;}

\itm[{\rm(ii)}] \textit{$\mbf{E}(m^{Y_0})\ge (m-1)\mbf{E}(Y_0m^{Y_0})$ implies $\sup_{n\ge 0} \mbf{E}(Y_n)\le (m-1)^{-1}$.}

 \eenumerate

According to the theorem above, the generalized DR process $(Y_n: n\ge 0)$ is \textit{supercritical}, \textit{subcritical} and \textit{critical} depending on the following conditions:
 \beqnn
\mbf{E}(m^{Y_0})~~ <,~ >~\mbox{and}~ = ~~ (m-1)\mbf{E}(Y_0m^{Y_0}).
 \eeqnn
The critical behavior of the free energy for the model was studied in Chen et al.~\cite{CDDHS21} and Hu and Shi~\cite{HuS18}, where weaker forms of \eqref{DR0conj} were established. The infinite smoothness of the free energy was proved by Chen~\cite{Che24+}.

Following Hu and Shi~\cite{HuS18}, we call $(Y_n: n\ge 0)$ a \textit{generalized DR process}. A review of the study of the model was given in Chen et al.~\cite{CDHLS19}, where the following conjectures were stated:

\medskip\noindent\textbf{Conjectures\,} {(Chen et al.~\cite{CDHLS19})\,} \textit{Suppose that $\eta\overset{\rm a.s.}= m\ge 2$. In the critical case, as $n\to \infty$,} \smallskip

\noindent{(C1)} ~\,\qqquad $\mbf{P}(Y_n\ge 1)= \frac{4}{(m-1)^2n^2} + o\big(\frac{1}{n^2}\big)$; \smallskip

\noindent{(C2)} ~\,\qqquad $\mbf{P}(Y_n= k|Y_n\geq 1)\to \frac{m-1}{m^k}$, $~k= 1, 2,\cdots$; \smallskip

\noindent{(C3)} ~\,\qqquad $\mbf{E}(Y_n)= \frac{4m}{(m-1)^3n^2} + o\big(\frac{1}{n^2}\big)$; \smallskip

\noindent{(C4a)} \qqquad $\mbf{E}(s^{Y_n})= 1 + \frac{4m(s-1)}{(m-1)^2(m-s)n^2} + o\big(\frac{1}{n^2}\big)$, $~0< s< m$; \smallskip

\noindent{(C4b)} \qqquad $\mbf{E}(m^{Y_n})= 1 + \frac{2}{(m-1)n} + o\big(\frac{1}{n}\big)$. \bigskip

The next theorem proved by Chen et al.~\cite{CHS22} establishes weaker forms of the conjectures~(C1) and~(C3) above.

\medskip\noindent\textbf{Theorem~B\,} {(Chen et al.~\cite{CHS22})\,} \textit{Suppose that $\eta\overset{\rm a.s.}=m\ge 2$ and $\mbf{E}[(m+\delta)^{Y_0}]< \infty$ for some $\delta> 0$. If $\mbf{E}(m^{Y_0})= (m-1)\mbf{E}(Y_0m^{Y_0})$, then, as $n\to \infty$,}
 \beqlb\label{P(Y_nge1)=1/n^2+o(.)}
\mbf{P}(Y_n\ge 1)= \frac{1}{n^{2+o(1)}},
 \quad
\mbf{E}(Y_n)= \frac{1}{n^{2+o(1)}}.
 \eeqlb

An exactly solvable continuous-time version of the DR model was introduced by Hu et al.~\cite{HMP20}, who gave complete characterizations for its phase transition and critical behavior. A scaling limit theorem for the generalized DR process was proved in Li and Zhang~\cite{LiZh24+}, which leads to a generalization of the model of Hu et al.~\cite{HMP20}.

In this work, we study the asymptotic behavior of the discrete-time generalized DR process with \textit{geometric offspring distribution}. More precisely, we assume the random variable $\eta$ in \eqref{gDR0proc} satisfies
 \beqlb\label{P(eta=n)=def}
\mbf{P}(\eta=n)= \frac{1}{m}\Big(1-\frac{1}{m}\Big)^{n-1}, \quad n=1,2,\cdots,
 \eeqlb
where $m= \mbf{E}(\eta)> 1$ is a constant. Given the parameters $(r,p)\in (0,1)^2 $, we denote by $\nu= G(r,p)$ the \textit{geometric-type distribution} given by:
 \beqnn
\nu(0)= p,~~ \nu(k)= (1-p)r(1-r)^{k-1}, \quad k= 1,2,\cdots.
 \eeqnn
For $n\ge 0$ let $\mu_n$ be the distribution of $Y_n$. Our main results are as follows.

\btheorem\label{th-mixgeo-mar0} Suppose that $\mu_0= G(r_0,p_0)$ for some $(r_0,p_0)\in (0,1)^2$. Then we have $\mu_n= G(r_n,p_n)$ for $n\geq 1$, where $\{(r_n,p_n): n\ge 1\}\subset (0,1)^2$ is defined recursively by
 \beqlb\label{mixgeo-mar1}
r_{n+1}= \frac{r_n}{m-(m-1)p_n},
\quad
p_{n+1}= 1 - (1-r_{n+1})\Big(1-\frac{r_{n+1}p_n}{r_n}\Big).
\eeqlb
\etheorem

\btheorem\label{th-(r_*,p_*)exist} Let $\{(r_n,p_n): n\ge 1\}\subset (0,1)^2$ be defined by \eqref{mixgeo-mar1} from any initial value $(r_0,p_0)\in (0,1)^2$. Then the following limits exist:
 \beqlb\label{r_*,p_*=def}
r_*:= \lim_{n\to \infty} r_n, \quad p_*:= \lim_{n\to \infty} p_n.
 \eeqlb
Moreover, one of the following holds:
 \beqnn
\text{$(1)$~$r_*=0$ and $p_*=0$; \quad $(2)$~$1-m^{-1}\le r_*< 1$ and $p_*=1$.}
 \eeqnn
\etheorem

\btheorem\label{th-mixgeoasy0} The following properties hold:
 \bitemize

\itm[{\rm(1)}] {\rm(Supercritical case)} When $r_*=0$ and $p_*=0$, we have
 \beqlb\label{F_infty}
F_\infty= \frac{1}{r_0}{\blue \prod_{i=0}^\infty }\Big[1-\frac{(m-1)p_i}{m}\Big]
 =
\frac{1-p_0}{r_0}\prod_{i=1}^{\infty} (1-r_i)\in (0,\infty)
 \eeqlb
and, as $n\to \infty$,
 \beqnn
\left\{\begin{aligned}
r_n&= \frac{1}{F_\infty m^n} - \frac{n}{F_\infty^2m^{2n-1}} + o\Big(\frac{n}{m^{2n}}\Big), \cr
p_n&= \frac{n}{F_\infty m^{n-1}} + \Big(\frac{p_0}{r_0} - mQ\Big)\frac{1}{F_\infty m^n} - \frac{n^2}{F_\infty^2m^{2(n-1)}} + o\Big(\frac{n^2}{m^{2n}}\Big),
\end{aligned}\right.
 \eeqnn
where $Q= \sum_{i=0}^\infty p_i< \infty$.

\itm[{\rm(2)}] {\rm(Subcritical case)} When $1-m^{-1}< r_*< 1$ and $p_*=1$, we have $F_\infty= 0$ and, as $n\to \infty$,
 \beqnn
\left\{\begin{aligned}
r_n&=  r_* + \frac{Kr_*^2\gamma_*^n}{1-\gamma_*} + \Big(1+\frac{mr_*}{1-\gamma_*^2}\Big) \frac{{\blue K^2}r_*^3\gamma_*^{2n}}{(1-\gamma_*)^2} + o(\gamma_*^{2n}), \cr
p_n&= 1 - \frac{Kr_*\gamma_*^n}{m-1} - {\blue\Big(1 +\frac{mr_*} {1-\gamma_*}\Big)} \frac{{\blue K^2}r_*^2\gamma_*^{2n}}{(m-1)(1-\gamma_*)} + o(\gamma_*^{2n}),
\end{aligned}\right.
 \eeqnn
where $\gamma_*= m(1-r_*)$ and
 \beqlb\label{K=..}
K= \frac{r_0-r_1}{r_0r_1}\prod_{i=1}^\infty \frac{1-r_i}{1-r_*}\in (0,\infty).
 \eeqlb

\itm[{\rm(3)}] {\rm(Critical case)} When $r_*= 1-m^{-1}$ and $p_*=1$, we have $F_\infty=0$ and, as $n\to \infty$,
 \beqnn
\left\{\begin{aligned}
r_n&= 1-\frac{1}{m} + \frac{2}{mn} - \frac{4(m+1)\log n}{3m(m-1)n^2} + o\Big(\frac{\log n}{n^2}\Big), \\
p_n&= 1 - \frac{2}{(m-1)^2n^2} - \frac{8(m+1)\log n}{3(m-1)^3 n^3} + o\Big(\frac{\log n}{n^3}\Big).
\end{aligned}\right.
 \eeqnn

 \eitemize
\etheorem

\bremark\label{re-iterativescheme0} By iterative arguments, we can extend the results of Theorem~\ref{th-mixgeoasy0} to higher-order asymptotic expansions. For example, the detailed arguments for the critical process are described in Remark~\ref{re-iterativescheme}. \eremark

As consequences of Theorem~\ref{th-mixgeoasy0}, we derive the asymptotics of the sustainability probability, first moment, marginal distribution and probability generating function as follows.

\bcorollary\label{th-sustprob0} As $n\to \infty$ we have:
 \bitemize

\itm[{\rm(1)}] in the supercritical case,
 \beqnn
\mbf{P}(Y_n\ge 1)= {\blue 1 - \frac{n}{F_\infty m^{n-1}} - \Big(\frac{p_0}{r_0} - mQ\Big)\frac{1}{F_\infty m^n} + \frac{n^2}{F_\infty^2m^{2(n-1)}} + o\Big(\frac{n^2}{m^{2n}}\Big);}
 \eeqnn

\itm[{\rm(2)}] in the subcritical case,
 \beqnn
\mbf{P}(Y_n\geq 1)\ar=\ar  {\blue \frac{Kr_*\gamma_*^n}{m-1} + \Big(1 -\frac{mr_*} {1-\gamma_*}\Big) \frac{ K^2 r_*^2\gamma_*^{2n}}{(m-1)(1-\gamma_*)} + o(\gamma_*^{2n});}
 \eeqnn

\itm[{\rm(3)}] in the critical case,
 \beqlb\label{P(Y_n>0)=asym}
\mbf{P}(Y_n\ge 1)= \frac{2}{(m-1)^2n^2} + \frac{8(m+1)\log n}{3(m-1)^3 n^3} + o\Big(\frac{\log n}{n^3}\Big).
 \eeqlb

 \eitemize
\ecorollary

\bcorollary \label{th-condlimitlaw0} As $n\to \infty$ we have:
\bitemize

\itm[{\rm(1)}] in the supercritical case, $Y_n\to \infty$ almost surely and
 \beqnn
\mbf{P}(m^{-n}Y_n\in \mathrm{d}x)\to F_{\infty}^{-1}\e^{-x F_{\infty}^{-1}}\mathrm{d}x, \quad x\ge 0;
 \eeqnn

\itm[{\rm(2)}] in the subcritical case, $Y_n\to 0$ almost surely and
 \beqnn
\mbf{P}(Y_n= k|Y_n\geq 1)\to (1-r_*)^{k-1}r_*, \quad k= 1, 2,\cdots;
 \eeqnn

\itm[{\rm(3)}] in the critical case, $Y_n\to 0$ almost surely and
 \beqlb\label{Y_infty=Geo}
\mbf{P}(Y_n= k|Y_n\geq 1)\to \Big(1-\frac{1}{m}\Big)\frac{1}{m^{k-1}}, \quad k= 1, 2,\cdots.
 \eeqlb

\eitemize
\ecorollary

\bremark The derivations of the results in Corollary~\ref{th-condlimitlaw0} essentially depends on the geometric structures of the offspring distribution and the marginal distributions. The exponential limit law in Corollary~\ref{th-condlimitlaw0}-(1) also appears in the scaling limit given by Hu et al.~\cite[Proposition~6.3]{HMP20} for their continuous-time model with exponential-type marginal distributions. Whether the exponential type limit holds for a generalized DR process with other offspring distributions or other marginal distributions still remains open. \eremark

\bcorollary\label{th-firmom0} As $n\to \infty$ we have:
 \bitemize

\itm[{\rm(1)}] in the supercritical case,
 \beqnn
\mbf{E}(Y_n)= F_\infty m^n + {\blue o(n)};
 \eeqnn

\itm[{\rm(2)}] in the subcritical case,
 \beqnn
\mbf{E}(Y_n)= \frac{K\gamma_*^n}{m-1} + {\blue \frac{mK^2r_*\gamma_*^{2n}}{(m-1)(1-\gamma_*)^2} } + o(\gamma_*^{2n});
 \eeqnn

\itm[{\rm(3)}] in the critical case,
 \beqlb\label{E(Y_n)=asym}
\mbf{E}(Y_n)= \frac{2m}{(m-1)^3n^2} + \frac{8m(m+1)\log n}{3(m-1)^4n^3} + o\Big(\frac{\log n}{n^3}\Big).
 \eeqlb

 \eitemize
\ecorollary

\bcorollary\label{th-E(.^{Y_n})-1=..0} As $n\to \infty$ we have:
 \bitemize

\itm[{\rm(1)}] in the supercritical case, for $|s|< 1$,
 \beqnn
\mbf{E}(s^{Y_n})= \frac{n}{F_\infty m^{n-1}} + \Big(\frac{p_0}{r_0} - mQ + \frac{s}{1-s}\Big)\frac{1}{F_\infty m^n} - \frac{n^2}{F_\infty^2m^{2(n-1)}} + o\Big(\frac{n^2}{m^{2n}}\Big);
 \eeqnn

\itm[{\rm(2)}] in the subcritical case, for $|s|< (1-r_*)^{-1}$,
 \beqnn
\mbf{E}(s^{Y_n})\ar=\ar{\blue  1 + \frac{s-1}{\,1-(1-r_*)s\,}
\Big\{
\frac{K r_*\,\gamma_*^n}{m-1}+\frac{ K^2\,r_*^2\gamma_*^{2n}}{(m-1)(1-\gamma_*)}}\ccr
\ar\ar{\blue\qquad \times
\left[
1+\frac{m r_*}{1-\gamma_*}
-\frac{sr_*}{\,1-(1-r_*)s\,}
\right]
+o(\gamma_*^{2n})
\Big\};}
\eeqnn

\itm[{\rm(3)}] in the critical case, for $|s|< m$,
 \beqlb\label{E(s^{Y_n})=1+..}
\mbf{E}(s^{Y_n})= 1 + \frac{s-1}{1-m^{-1}s} \Big[\frac{2}{(m-1)^2n^2} + \frac{8(m+1)\log n}{3(m-1)^3 n^3} + o\Big(\frac{\log n}{n^3}\Big)\Big],
 \eeqlb
and
 \beqlb\label{E(m^{Y_n})=1+..}
\mbf{E}(m^{Y_n})= 1 + \frac{1}{(m-1)n} + \frac{2(m+1)\log n}{(m-1)^2n^2} + o\Big(\frac{\log n}{n^2}\Big).
 \eeqlb

 \eitemize
\ecorollary

\bremark\label{th-poss-improvs1} In the proofs of the above corollaries, we only use the exact forms of the first infinitesimal terms of the expansions Theorem~\ref{th-mixgeoasy0}. The results can be improved by making full use of the expansions. \eremark

\bremark\label{th-dif-conj-result0} The distribution in \eqref{Y_infty=Geo} agrees with that in (C2), but the coefficients of the first order infinitesimal terms in \eqref{P(Y_n>0)=asym}, \eqref{E(Y_n)=asym}, \eqref{E(s^{Y_n})=1+..} and~\eqref{E(m^{Y_n})=1+..} do not agree exactly with those in (C1), (C3), (C4a) and~(C4b). This seems due to the random fluctuations of the offspring number $\eta$, which slow down the mass production compared to the deterministic offspring number case of the original conjectures. \eremark

\bremark\label{th-mixexp-mar0} Given the parameters $(\lambda,p)\in (0,\infty)\times (0,1)$, we denote by $\mu= E(\lambda,p)$ the \textit{exponential-type distribution} given by:
 \beqnn
\mu(\d x)= p\delta_0+(1-p)\lambda\e^{-\lambda x}\d x, \quad x\ge 0,
 \eeqnn
where $\delta_0$ is the unit mass at zero. As in the proof of Theorem~\ref{th-mixgeo-mar0}, one can show that if $\mu_0= E(\lambda_0,p_0)$ for $(\lambda_0,p_0)\in (0,\infty)\times (0,1)$, then $\mu_n= E(\lambda_n,p_n)$ for $n\geq 1$, where the sequence $\{(\lambda_n,p_n): n\ge 1\}\subset (0,\infty)\times (0,1)$ is defined recursively by
 \beqnn
\lambda_{n+1}= \frac{\e^{-\alpha}\lambda_n}{1-(1-\e^{-\alpha})p_n},
 \quad
p_{n+1}= 1 - \e^{-\lambda_{n+1}}\Big(1 - \frac{\lambda_{n+1}p_n}{\lambda_n}\Big).
 \eeqnn
All the results given above can be extended to the exponential-type marginal distributions by similar arguments. \eremark

\bremark\label{th-regimes0} The results of Theorem~\ref{th-mixgeoasy0} are not quite satisfactory as the supercritical, subcritical and critical regimes are characterized by $(r_*,p_*)$, rather than $(r_0,p_0)$. By \eqref{mixgeo-mar1} the sequence $\{r_n\}\subset (0,1)$ is strictly decreasing. Then we have $\lim_{n\to \infty} r_n= 0$ when $0< r_0\le 1-m^{-1}$. This means that the generalized DR model belongs to the supercritical regime if $(r_0,p_0)\in D:= (0,1-m^{-1}]\times (0,1)$. A further generalization of the recursive model \eqref{gDR0proc} was studied in the very recent work of Alsmeyer et al.~\cite{AHM25+}. Their model still preserves the geometric and exponential-type marginal distributions. For those marginal distributions, they proved the Derrida--Retaux conjecture \eqref{DR0conj} and gave a very ingenious characterization of the critical curve using an involution-type equation. \eremark

The rest of the paper is organized as follows. The basic structures of the geometric-type marginal distributions are discussed in Section~2, where the proofs of Theorems~\ref{th-mixgeo-mar0} and~\ref{th-(r_*,p_*)exist} are also given. The proofs of Theorem~\ref{th-mixgeoasy0} and its corollaries are given in Section~3.

\section{Geometric-type marginal distributions}

\setcounter{equation}{0}

\medskip\noindent\textit{Proof of Theorem~\ref{th-mixgeo-mar0}.~} We show the result by induction in $n\geq 0$. By the assumption, we have $\mu_0= G(r_0,p_0)$. Now suppose that $\mu_n= G(r_n,p_n)$ for some $n\geq 0$. Then, for $0< s< (1-r_n)^{-1}$,
 \beqlb\label{E(s^{Y_n})=p_n+..}
\mbf{E}(s^{Y_n})
 =
p_n + \frac{(1-p_n)r_ns}{1-(1-r_n)s}
 =
\frac{p_n + (r_n-p_n)s}{1-(1-r_n)s}.
 \eeqlb
Write $\xi_{n+1}= \sum_{k=1}^{\eta_n} Y_{n,k}$. By the independence of $\eta_n$ and $\{Y_{n,k}\}$, we see that
 \beqnn
\mbf{E}(s^{-\xi_{n+1}})\ar=\ar \sum_{k=1}^\infty m^{-1}(1-m^{-1})^{k-1} \Big[\frac{p_n + (r_n-p_n)s}{1-(1-r_n)s}\Big]^k \cr
 \ar=\ar
\frac{p_n + (r_n-p_n)s}{m-(m-1)p_n-[m - (m-1)p_n - r_n]s} \ccr
 \ar=\ar
q_{n+1} + (1-q_{n+1})\frac{r_{n+1} s}{1-(1-r_{n+1})s},
 \eeqnn
where
 \beqnn
r_{n+1}= \frac{r_n}{m-(m-1)p_n},
 \quad
q_{n+1}= \frac{p_n}{m-(m-1)p_n}= \frac{r_{n+1}p_n}{r_n}.
 \eeqnn
Then $\xi_{n+1}$ follows the geometric-type distribution $G(r_{n+1},q_{n+1})$, that is,
 \beqnn
\mbf{P}(\xi_{n+1}\in\d x)= q_{n+1}\delta_0(\d x) + (1-q_{n+1})r_{n+1}\sum_{k=1}^\infty(1-r_{n+1})^{k-1}\delta_k(\d x).
 \eeqnn
By the total probability formula and the memoryless of the geometric distribution, we have
 \beqnn
\mbf{P}(Y_{n+1}\in\d x)\ar=\ar \mbf{P}(\xi_{n+1}\ge 2)\mbf{P}(\xi_{n+1}-1\in\d x|\xi_{n+1}\ge 2) + \mbf{P}(\xi_{n+1}\le 1)\delta_0(\d x) \ccr
 \ar=\ar
(1-q_{n+1})(1-r_{n+1})\sum_{k=1}^\infty r_{n+1}(1-r_{n+1})^{k-1}\delta_k(\d x) \ccr
 \ar\ar
+\, [1-(1-r_{n+1})(1-q_{n+1})]\delta_0(\d x),
 \eeqnn
where
 \beqnn
(1-q_{n+1})(1-r_{n+1})= (1-r_{n+1})\Big(1 - \frac{r_{n+1}p_n}{r_n}\Big)= 1-p_{n+1}.
 \eeqnn
Then we have $\mu_{n+1}= G(r_{n+1},p_{n+1})$. \qed

\bproposition\label{th-m^{-n}E(Y_n)(1)} Let $\{(r_n,p_n): n\ge 1\}$ be the sequence defined by \eqref{mixgeo-mar1}. Then for any $n\geq 0$ we have
 \beqlb\label{m^{-n}E(Y_n)=..}
m^{-n}\mbf{E}(Y_n)= r_0^{-1}(1-p_n){\blue \prod_{i=0}^{n-1}} \Big[1-\frac{(m-1)p_i}{m}\Big]
 =
r_0^{-1}(1-p_0) \prod_{i=1}^n (1-r_i).
 \eeqlb
\eproposition

\bproof By Theorem~\ref{th-mixgeo-mar0} we know that $Y_n$ has the geometric-type distribution $G(r_n,p_n)$. It follows that
 \beqlb\label{m^{-n}E(Y_n)=(0)}
\mbf{E}(Y_n)= r_n^{-1}(1-p_n).
 \eeqlb
By the first equality in \eqref{mixgeo-mar1}, we have
 \beqlb\label{r_n/r_{n+1}=..}
\frac{r_n}{r_{n+1}}= m-(m-1)p_n,
\eeqlb
 and hence
 \beqlb\label{r_n=r_0prod}
r_n= r_0{\blue\prod_{i=0}^{n-1}} \big[m-(m-1)p_i\big]^{-1}.
 \eeqlb
Then the first equality in \eqref{m^{-n}E(Y_n)=..} follows by \eqref{m^{-n}E(Y_n)=(0)}. By the two equalities in \eqref{mixgeo-mar1} we have
 \beqlb\label{p_{n+1}=1-..(2)}
p_{n+1}\ar=\ar 1 - (1-r_{n+1})\Big(1 - \frac{r_{n+1}p_n}{r_n}\Big) \cr
 \ar=\ar
1 - (1-r_{n+1})\Big(1 - \frac{p_n}{m-(m-1)p_n}\Big) \cr
 \ar=\ar
1 - \frac{m(1-p_n)(1-r_{n+1})}{m-(m-1)p_n}.
 \eeqlb
It follows that
 \beqlb\label{1-p_{n+1}=..}
1-p_{n+1}= m(1-p_n)\frac{r_{n+1}}{r_n}(1-r_{n+1}).
 \eeqlb
and hence
 \beqnn
\frac{1-p_{n+1}}{r_{n+1}}= m(1-r_{n+1})\frac{1-p_n}{r_n}.
 \eeqnn
This proves the second equality in \eqref{m^{-n}E(Y_n)=..}. \eproof

\bproposition\label{th-eq-r_n} Let $\{(r_n,p_n): n\ge 1\}$ be the sequence defined by \eqref{mixgeo-mar1}. Then for any $n\geq 0$ we have
 \beqlb\label{eq-r_n}
\Big(\frac{1}{r_{n+2}}-\frac{1}{r_{n+1}}\Big)
 =
m(1-r_{n+1})\Big(\frac{1}{r_{n+1}}-\frac{1}{r_n}\Big).
 \eeqlb
\eproposition

\bproof From \eqref{r_n/r_{n+1}=..} it follows that
 \beqlb\label{1-p_n=..}
1-p_n=\frac{1}{m-1}\Big(\frac{r_n}{r_{n+1}}-1\Big).
 \eeqlb
Substituting the above expression into \eqref{1-p_{n+1}=..} we obtain
 \beqnn
\frac{1}{m-1}\Big(\frac{r_{n+1}}{r_{n+2}}-1\Big)
 \ar=\ar
\frac{m}{m-1}\Big(\frac{r_n}{r_{n+1}}-1\Big)\frac{r_{n+1}}{r_n}(1-r_{n+1}) \cr
 \ar=\ar
\frac{m}{m-1}\Big(1-\frac{r_{n+1}}{r_n}\Big)(1-r_{n+1}).
 \eeqnn
This proves \eqref{eq-r_n}. \eproof

\bproposition\label{th-eq-r_np_n} Let $\{(r_n,p_n): n\ge 1\}$ be the sequence defined by \eqref{mixgeo-mar1}. Then for any $n\geq 0$ we have
 \beqlb\label{eq-r_np_n}
\frac{p_{n+1}}{r_{n+1}}-\frac{p_n}{r_n}= m(1-p_n).
 \eeqlb
\eproposition

\bproof By \eqref{r_n/r_{n+1}=..} and \eqref{p_{n+1}=1-..(2)} it follows that
 \beqnn
p_{n+1}\ar=\ar \frac{m(1-p_n)+p_n-m(1-r_{n+1})(1-p_n)}{m-(m-1)p_n} \cr
 \ar=\ar
\frac{m(1-p_n)r_n+p_n}{m-(m-1)p_n}
 =
\frac{r_{n+1}}{r_n}[m(1-p_n)r_n+p_n] \cr
 \ar=\ar
r_{n+1}\Big[m(1-p_n) + \frac{p_n}{r_n}\Big].
 \eeqnn
Then the desired relation holds. \eproof

\medskip\noindent\textit{Proof of Theorem~\ref{th-(r_*,p_*)exist}.~} By \eqref{mixgeo-mar1} it is easy to see the sequence $\{r_n\}$ is strictly decreasing. Set $r_*= \lim_{n\to \infty} r_n$. We claim that either $r_* = 0$ or $r_*\in [1-m^{-1},\infty)$. Otherwise, we have $r_*\in (0,1-m^{-1})$, and so there exists $N\ge 1$ such that $0< r_*\le r_n< 1-m^{-1}$ for all $n\ge N$. In this case, it follows from \eqref{eq-r_n} that
 \beqnn
\Big|\frac{1}{r_{n+2}}-\frac{1}{r_{n+1}}\Big|> \Big|\frac{1}{r_{n+1}}-\frac{1}{r_n}\Big|,
 \eeqnn
which is in contradiction to the existence of the limit $\lim_{n\to \infty} 1/r_n= 1/r_*$. Therefore, we must have $r_* = 0$ or $r_*\in [1-m^{-1},\infty)$. We next examine the convergence of the sequence $\{p_n\}$ in the two cases.

(1)~In the case $r_*= 0$, the sequence $1/r_n$ strictly increases to $\infty$ as $n\to \infty$. By \eqref{eq-r_n} and the Stolz--Ces\`aro theorem, we have
 \beqlb\label{r_n/r_{n+1}to m}
\lim_{n\to \infty} \frac{r_{n+1}}{r_n}
 =
\lim_{n\to \infty} \frac{\frac{1}{r_n}}{\frac{1}{r_{n+1}}}
 =
\lim_{n\to \infty} \frac{\frac{1}{r_n}-\frac{1}{r_{n-1}}} {\frac{1}{r_{n+1}}-\frac{1}{r_n}}
 =
\lim_{n\to \infty}m^{-1}(1-r_n)^{-1}= m^{-1}.
 \eeqlb
Then $p_*= \lim_{n\to \infty} p_n= 0$ by \eqref{r_n/r_{n+1}=..}.

(2)~In the case $r_*\in [1-m^{-1},\infty)$, we have $\lim_{n\to \infty} (r_n/r_{n+1})= 1$ and hence $p_*= \lim_{n\to \infty} p_n= 1$ by \eqref{1-p_n=..}. \qed

\section{Asymptotic behavior of the dynamics}

\setcounter{equation}{0}

\medskip\noindent\textit{Proof of Theorem~\ref{th-mixgeoasy0}-(1).~} In the supercritical case, the limits defined in \eqref{r_*,p_*=def} satisfy $r_*= 0$ and $p_*= 0$. From \eqref{r_n/r_{n+1}to m} we see that
 \beqnn
0 < \prod_{n=1}^\infty (1-r_n)< \infty.
 \eeqnn
Then we get \eqref{F_infty} from \eqref{m^{-n}E(Y_n)=..}. From \eqref{eq-r_n} it follows that
 \beqlb\label{(1/r_{n+1}-1/r_n)=..}
\frac{1}{r_{n+1}}-\frac{1}{r_n}
 =
m^n \Big(\frac{1}{r_{1}}-\frac{1}{r_{0}}\Big)\prod_{i=1}^n (1-r_i),
 \eeqlb
Since $r_n^{-1}$ increases strictly to $\infty$ as $n\to\infty$, in view of \eqref{r_n/r_{n+1}to m} and \eqref{(1/r_{n+1}-1/r_n)=..}, we can use the Stolz--Ces\`aro theorem to see that
 \beqnn
\lim_{n\to\infty} \frac{r_n}{m^{- n}}
 \ar=\ar
\lim_{n\to\infty} \frac{m^n}{\frac{1}{r_n}}
=
\lim_{n\to\infty} \frac{m^n(m-1)}{\frac{1}{r_{n+1}}-\frac{1}{r_n}} \cr
 \ar=\ar
 \frac{m-1}{\big(\frac{1}{r_{1}} - \frac{1}{r_0}\big)}\prod_{n=1}^\infty (1-r_n)^{-1}
 =
F_\infty^{-1},
 \eeqnn
which implies
 \beqlb\label{r_n=..}
r_n= F_\infty^{-1}m^{-n}+o(m^{-n}).
 \eeqlb
By using \eqref{eq-r_np_n}, \eqref{r_n/r_{n+1}to m} and the Stolz--Ces\`aro theorem, we have
 \beqnn
\lim_{n\to\infty}\frac{p_n}{nr_n}\ar=\ar \lim_{n\to\infty}\frac{\frac{p_n}{r_n}}{n}
 =
\lim_{n\to\infty}\bigg(\frac{p_{n+1}}{r_{n+1}}-\frac{p_n}{r_n}\bigg)
 =
\lim_{n\to\infty}m(1-p_n)= m.
 \eeqnn
This together with \eqref{r_n=..} gives
 \beqlb\label{p_n=..}
p_n= F_\infty^{-1}nm^{1-n} + o(nm^{-n}).
 \eeqlb
We next deduce the high-order expansion of $(r_n,p_n)$. By the Stolz--Ces\`aro theorem and \eqref{p_n=..}, it is simple to see that, as $n\to \infty$,
 \beqlb\label{sump_i=..}
\sum_{i=n}^\infty p_i
 \ar=\ar
[1+o(1)]\sum_{i=n}^\infty \frac{im^{1-i}}{F_\infty}
 =
[1+o(1)]\frac{m^{2-n}(nm-n+1)}{F_\infty(m-1)^2} \cr
 \ar=\ar
\frac{nm^{2-n}}{F_\infty(m-1)} + o\Big(\frac{n}{m^n}\Big),
 \eeqlb
and hence,
 \beqlb\label{sumlog[1-..]=..}
\sum_{i=n}^\infty \log\Big[1-\frac{(m-1)p_i}{m}\Big]
 \ar=\ar
-[1+o(1)]\sum_{i=n}^\infty \frac{(m-1)p_i}{m} \cr
 \ar=\ar
-\frac{n}{F_\infty m^{n-1}} + o\Big(\frac{n}{m^n}\Big).
 \eeqlb
By \eqref{eq-r_np_n} we have
 \beqnn
\frac{p_n}{r_n}= \frac{p_0}{r_0} + m\sum_{i=0}^{n-1} (1-p_i)
 =
\frac{p_0}{r_0} + mn - mQ + m\sum_{i=n}^\infty p_i.
 \eeqnn
It follows that
 \beqlb\label{p_n=..(2)}
p_n= \Big[mn + \Big(\frac{p_0}{r_0} - mQ\Big) + m\sum_{i=n}^\infty p_i\Big]r_n.
 \eeqlb
From \eqref{F_infty} and \eqref{r_n=r_0prod} one can see that
 \beqnn
r_n\ar=\ar \frac{r_0}{m^n} \prod_{i=0}^{n-1} \Big[1-\frac{(m-1)p_i}{m}\Big]^{-1} \cr
 \ar=\ar
\frac{1}{F_\infty m^n} \exp\Big\{\sum_{i=n}^\infty \log\Big[1-\frac{(m-1)p_i}{m}\Big]\Big\} \cr
 \ar=\ar
\frac{1}{F_\infty m^n} \exp\Big\{-\frac{n}{F_\infty m^{n-1}} + o\Big(\frac{n}{m^n}\Big)\Big\} \cr
 \ar=\ar
\frac{1}{F_\infty m^n} \Big[1 - \frac{n}{F_\infty m^{n-1}} + o\Big(\frac{n}{m^n}\Big)\Big],
 \eeqnn
where we used \eqref{sumlog[1-..]=..} and Taylor's expansion of the exponential function. This gives the desired expansion of $r_n$. Furthermore, by \eqref{sump_i=..} and \eqref{p_n=..(2)} we have
 \beqnn
p_n\ar=\ar \Big[mn + \Big(\frac{p_0}{r_0} - mQ\Big) + m\sum_{i=n}^\infty p_i\Big] \cdot \Big[\frac{1}{F_\infty m^n} - \frac{n}{F_\infty^2m^{2n-1}} + o\Big(\frac{n}{m^{2n}}\Big)\Big] \cr
 \ar=\ar
\frac{n}{F_\infty m^{n-1}} + \Big(\frac{p_0}{r_0} - mQ\Big)\frac{1}{F_\infty m^n} - \frac{n^2}{F_\infty^2m^{2(n-1)}} + o\Big(\frac{n^2}{m^{2n}}\Big),
 \eeqnn
which gives the desired expansion of $p_n$. \qed

\medskip\noindent\textit{Proof of Theorem~\ref{th-mixgeoasy0}-(2).~} In the subcritical case, the limits defined in \eqref{r_*,p_*=def} satisfy $1-m^{-1}< r_*< 1$ and $p_*= 1$.  By \eqref{m^{-n}E(Y_n)=..} we have $F_\infty= \lim_{n\to \infty} m^{-n}\mbf{E}(Y_n)= 0$. From \eqref{eq-r_n} it follows that
 \beqlb\label{1/lam_{n+2}-1/lam_{n+1}}
r_n-r_{n+1}= m(1-r_n) \frac{r_{n+1}}{r_{n-1}} (r_{n-1}-r_n).
 \eeqlb
Then an application of the Stolz--Ces\`aro theorem leads to
 \beqnn
\lim_{n\to \infty} \frac{r_{n+1}-r_*}{r_n-r_*}
 \ar=\ar
\lim_{n\to \infty} \frac{r_{n+1}-r_n}{r_n-r_{n-1}} \cr
 \ar=\ar
\lim_{n\to \infty} m(1-r_n) \frac{r_{n+1}}{r_{n-1}}= \gamma_*< 1,
 \eeqnn
which implies that
 \beqnn
0< \prod_{i=1}^\infty \frac{1-r_i}{1-r_*}< \infty.
 \eeqnn
By using \eqref{1/lam_{n+2}-1/lam_{n+1}} again we deduce that
 \beqlb\label{r_n-r_{n+1}=m^n..}
r_n-r_{n+1}\ar=\ar m^n\frac{r_nr_{n+1}}{r_0r_1}(r_0-r_1)\prod_{i=1}^n (1-r_i) = K_n \gamma_*^n r_nr_{n+1},
 \eeqlb
where
 \beqlb\label{K_n=..}
K_n= \frac{r_0-r_1}{r_0r_1}\prod_{i=1}^n \frac{1-r_i}{1-r_*}.
 \eeqlb
It follows that
 \beqnn
\lim_{n\to \infty} \frac{r_n-r_{n+1}}{\gamma_*^n}= Kr_*^2,
 \eeqnn
where $K$ is defined by \eqref{K=..}. By the Stolz--Ces\`aro theorem,
 \beqnn
\lim_{n\to \infty} \frac{r_n-r_*}{\gamma_*^n }=\lim_{n\to \infty} \frac{r_n-r_{n+1}}{(1-\gamma_*)\gamma_*^n} = \frac{Kr_*^2}{1-\gamma_*},
 \eeqnn
which implies that
 \beqlb\label{r_n-r_*=gamma_*^n...}
r_n= r_* + \frac{Kr_*^2\gamma_*^n}{1-\gamma_*} + o(\gamma_*^n).
 \eeqlb
By the Stolz--Ces\`aro theorem, as $n\to \infty$,
 \beqnn
\sum_{i=n+1}^\infty \log \Big(1-\frac{r_i-r_*}{1-r_*}\Big)
 \ar=\ar
-[1+o(1)]\sum_{i=n+1}^\infty \frac{r_i-r_*}{1-r_*} \cr
 \ar=\ar
-[1+o(1)]\sum_{i=n+1}^\infty \frac{}{}\frac{Kr_*^2\gamma_*^{\blue i}}{(1-r_*)(1-\gamma_*)} \cr
 \ar=\ar
-\frac{mKr_*^2\gamma_*^n}{(1-\gamma_*)^2} + o(\gamma_*^n).
 \eeqnn
Then by \eqref{K=..} and \eqref{K_n=..} we can write
 \beqnn
K_n\ar=\ar K\exp\Big\{-\sum_{i=n+1}^\infty \log\Big(1-\frac{r_i -r_*}{1-r_*}\Big)\Big\} \cr
 \ar=\ar
K\exp\Big\{\frac{mKr_*^2\gamma_*^n}{(1-\gamma_*)^2} + o(\gamma_*^n)\Big\} \cr
 \ar=\ar
K + \frac{mK^2r_*^2\gamma_*^n}{(1-\gamma_*)^2} + o(\gamma_*^n),
 \eeqnn
where the last equality holds by the Taylor's expansion of the exponential function. Substituting this and \eqref{r_n-r_*=gamma_*^n...} into \eqref{r_n-r_{n+1}=m^n..} we get
 \beqlb\label{r_n-r_{n+1}=m^n...}
r_n-r_{n+1}\ar=\ar \gamma_*^n\Big[K + \frac{mK^2r_*^2\gamma_*^n}{(1-\gamma_*)^2} + o(\gamma_*^n)\Big]\cdot \Big[r_* + \frac{Kr_*^2\gamma_*^n}{1-\gamma_*} + o(\gamma_*^n)\Big] \cr
 \ar\ar\quad
\cdot\,\Big[r_* + \frac{Kr_*^2\gamma_*^{n+1}}{1-\gamma_*} + o(\gamma_*^n)\Big] \cr
 \ar=\ar
Kr_*^2\gamma_*^n + \Big(1+\gamma_* + \frac{mr_*}{1-\gamma_*}\Big) \frac{{\blue K^2}r_*^3\gamma_*^{2n}}{1-\gamma_*} + o(\gamma_*^{2n}).
 \eeqlb
By the Stolz--Ces\`aro theorem, as $n\to \infty$ we have
 \beqlb\label{r_n-r_*=...}
r_n\ar=\ar r_* + \sum_{i=n}^\infty (r_i-r_{i+1}) \cr
 \ar=\ar
r_* + [1+o(1)]\sum_{i=n}^\infty \Big[Kr_*^2\gamma_*^{\blue i} + \Big(1+\gamma_* + \frac{mr_*}{1-\gamma_*}\Big) \frac{{\blue K^2}r_*^3\gamma_*^{{\blue 2i}}}{1-\gamma_*}\Big] \cr
 \ar=\ar
r_* + \frac{Kr_*^2\gamma_*^n}{1-\gamma_*} + \Big(1 + \frac{mr_*}{1-\gamma_*^2}\Big) \frac{{\blue K^2}r_*^3\gamma_*^{2n}}{(1-\gamma_*)^2} + o(\gamma_*^{2n}),
 \eeqlb
which gives the desired expansion of $r_n$. By \eqref{r_n-r_*=...} and Taylor's expansion,
 \beqnn
\frac{1}{r_n}= \frac{1}{r_* + (r_n-r_*)}
 =
\frac{1}{r_*} - \frac{r_n-r_*}{r_*^2} + o(r_n-r_*)
 =
\frac{1}{r_*} - \frac{K\gamma_*^n}{1-\gamma_*} + o(\gamma_*^n).
 \eeqnn
Using \eqref{1-p_n=..}, \eqref{r_n-r_{n+1}=m^n...} and the above expansion, we obtain
 \beqnn
1-p_n\ar=\ar \frac{1}{m-1}(r_n-r_{n+1}) \frac{1}{r_{n+1}} \cr
 \ar=\ar
\frac{1}{m-1}\Big[Kr_*^2\gamma_*^n + \Big(1+\gamma_* + \frac{mr_*}{1-\gamma_*}\Big) \frac{{\blue K^2}r_*^3\gamma_*^{2n}}{1-\gamma_*} + o(\gamma_*^{2n})\Big] \cr
 \ar\ar\qquad
\cdot\, \Big[\frac{1}{r_*} - \frac{K\gamma_*^{n+1}}{1-\gamma_*} + o(\gamma_*^n)\Big] \cr
 \ar=\ar
\frac{Kr_*\gamma_*^n}{m-1} + {\blue\Big(1 +\frac{mr_*} {1-\gamma_*}\Big)} \frac{{\blue K^2}r_*^2\gamma_*^{2n}}{(m-1)(1-\gamma_*)} + o(\gamma_*^{2n}).
 \eeqnn
Then we have the desired expansion of $p_n$. \qed

We need some preparations for the proof of Theorem~\ref{th-mixgeoasy0}-(3), which deals with the critical process. In this case, the limits defined in \eqref{r_*,p_*=def} satisfy $r_*= 1-m^{-1}$ and $p_*= 1$. To simplify the presentation, let us consider the sequence
 \beqa\label{eq:v_n=r_n-(1-m^{-1})}
v_n = r_n - (1 - m^{-1}), \quad n\ge 0.
 \eeqa

\blemma\label{th(nv_n)(n^2(v_n-v_{n+1}))} In the critical case, we have
 \beqlb\label{(nv_n)(n^2(v_n-v_{n+1}))}
\lim_{n\to \infty} nv_n= \frac{2}{m},
 \quad
\lim_{n\to \infty} n^2(v_n-v_{n+1})
 =
\lim_{n\to \infty} \frac{m}{2}n^2v_nv_{n+1}= \frac{2}{m}.
 \eeqlb
\elemma

\bproof To simplify the presentation, we introduce the difference operator $\Delta$ in the following way: {\blue for} any sequence $\{a_n\}$ write $\Delta a_n= a_{n+1}-a_n$ and
 \beqnn
\Delta^2 a_n=\Delta(\Delta a_n)= a_{n+2}-2a_{n+1}+a_n.
 \eeqnn
By \eqref{eq-r_n} we have
 \beqnn
(r_{n+2}-r_{n+1})
 \ar=\ar
m(1-r_{n+1})\frac{r_{n+2}}{r_n}(r_{n+1}-r_n) \cr
 \ar=\ar
(1-m v_{n+1})\frac{r_{n+2}}{r_n}(r_{n+1}-r_n),
 \eeqnn
which implies
 \beqlb\label{eqDeltv_{n+1}}
\Delta v_{n+1}= (1-mv_{n+1}) \frac{v_{n+2}+1-m^{-1}}{v_n+1-m^{-1}}\Delta v_n.
\eeqlb
It follows that
 \beqlb\label{eqDelt^2v_n}
\Delta^2 v_n\ar=\ar (1-mv_{n+1}) \frac{v_{n+2}+1-m^{-1}}{v_n+1-m^{-1}}\Delta v_n - \Delta v_n \cr
 \ar=\ar
\Big[-m v_{n+1} + (1-m v_{n+1})\frac{v_{n+2}-v_n}{v_n+1-m^{-1}}\Big] \Delta v_n.
 \eeqlb
Note that $v_n$ strictly decreases to zero as $n\to \infty$. By applying the Stolz--Ces\`aro theorem,
 \beqnn
\lim_{n\to \infty} \frac{v_{n+1}}{v_n}
 =
\lim_{n\to \infty} \frac{\Delta v_{n+1}}{\Delta v_n}
 =
\lim_{n\to \infty} (1-m v_{n+1}) \frac{v_{n+2}+1-m^{-1}}{v_n+1-m^{-1}}= 1.
 \eeqnn
Then we deduce
 \beqlb\label{limfrac{1}{v_n}[..]}
\ar\ar\lim_{n\to \infty} \frac{1}{v_n} \Big[(1-m v_{n+1}) \Big(1 + \frac{v_{n+2}-v_n}{v_n+1-m^{-1}}\Big) - 1\Big] \cr
 \ar\ar\qquad
= \lim_{n\to \infty} \Big(-m\frac{v_{n+1}}{v_n}+ \frac{v_{n+2}/v_n-1}{v_n+1-m^{-1}}- \frac{v_n-v_{n+2}}{v_n+1-m^{-1}}\frac{v_{n+1}}{v_n}\Big) = -m.
 \eeqlb
By \eqref{eqDeltv_{n+1}} the sequence $-\Delta v_{n+1}$ strictly decreases to zero. Then we can use \eqref{limfrac{1}{v_n}[..]} and the Stolz--Ces\`aro theorem to obtain
 \beqnn
\lim_{n\to \infty}\frac{v_{n+1}v_n}{v_n-v_{n+1}}
 \ar=\ar
\lim_{n\to \infty}\frac{v_{n+1}(v_{n+2}-v_n)}{(v_{n+1}-v_{n+2})-(v_n-v_{n+1})} \cr
 \ar=\ar
\lim_{n\to \infty}\frac{v_{n+1}(\Delta v_n+\Delta v_{n+1})}{-\Delta^2 v_n}\cr
 \ar=\ar
\lim_{n\to \infty}\frac{-v_n(1+\Delta v_{n+1}/\Delta v_n)}{(1-m v_{n+1}) \big(1+\frac{v_{n+2}-v_n}{v_n+1-m^{-1}}\big)-1}= \frac{2}{m}.
 \eeqnn
By another application of the Stolz--Ces\`aro theorem we get
 \beqnn
\lim_{n\to \infty} nv_n
 =
\lim_{n\to \infty} \frac{1}{\frac{1}{v_{n+1}}-\frac{1}{v_n}}
 =
\lim_{n\to \infty} \frac{v_nv_{n+1}}{v_n-v_{n+1}}= \frac{2}{m}.
 \eeqnn
This gives the limits in \eqref{(nv_n)(n^2(v_n-v_{n+1}))}. \eproof

From the results of Lemma~\ref{th(nv_n)(n^2(v_n-v_{n+1}))}, we can derive a refined asymptotic expansion for $v_n$ as $n\to \infty$. The the approach depends on the basic relation:
 \beqa\label{eq:1/v_{n+1}-1/v_n-m/2=...}
\Big(\frac{1}{v_{n+1}}-\frac{1}{v_n}\Big) - \frac{m}{2}
 =
\frac{1}{v_nv_{n+1}}\Big(v_n - v_{n+1} - \frac{m}{2}v_nv_{n+1}\Big).
 \eeqa

\blemma\label{thv_n-v_{n+1}} In the critical case, we have
 \beqlb\label{eq:v_n-v_{n+1}}
\lim_{n\to \infty} n^3\Big(v_n - v_{n+1} - \frac{m}{2}v_nv_{n+1}\Big)
 =
\frac{4(m+1)}{3m(m-1)}.
 \eeqlb
\elemma

\bproof Let the difference operator $\Delta$ be defined in the proof of Lemma~\ref{th(nv_n)(n^2(v_n-v_{n+1}))}. In view of \eqref{(nv_n)(n^2(v_n-v_{n+1}))}, as $n\to \infty$ we have
 \beqnn
v_n= \frac{2}{mn} + o\Big(\frac{1}{n}\Big),
 \quad
\Delta v_n= -\frac{2}{mn^2} + o\Big(\frac{1}{n^2}\Big).
 \eeqnn
Then by Taylor's expansion for the function $1/(1+x)$ we see that
 \beqnn
\frac{v_{n+2}-v_n}{v_n+1-m^{-1}}\ar=\ar \frac{m}{m-1}(\Delta v_n+\Delta v_{n+1})\frac{1}{1+mv_n/(m-1)} \cr
 \ar=\ar
\frac{m}{m-1}(\Delta v_n+\Delta v_{n+1})\Big[1 - \frac{m}{m-1}v_n + o\Big(\frac{1}{n}\Big)\Big] \cr
 \ar=\ar
\frac{m}{m-1}(\Delta v_n+\Delta v_{n+1}) + o\Big(\frac{1}{n^2}\Big).
 \eeqnn
Substituting the above expression into \eqref{eqDelt^2v_n} we obtain
 \beqnn
\Delta^2 v_n= -mv_{n+1}\Delta v_n + \frac{m}{m-1}(\Delta v_n+\Delta v_{n+1})\Delta v_n + o\Big(\frac{1}{n^4}\Big).
 \eeqnn
It follows that
 \beqlb\label{eq:Delta ( Delta v_n+frac{m}{2}v_{n+1}v_n )=...}
\Delta \Big( \Delta v_n+\frac{m}{2}v_{n+1}v_n\Big)
 \ar=\ar
\Delta^2 v_n + \frac{m}{2}v_{n+1}\Delta v_{n+1} + \frac{m}{2}v_{n+1}\Delta v_n \cr
 \ar=\ar
- \frac{m}{2}v_{n+1}\Delta v_n + \frac{m}{m-1}(\Delta v_n+\Delta v_{n+1})\Delta v_n \cr
 \ar\ar
+\,\frac{m}{2}v_{n+1}\Delta v_{n+1} + o\Big(\frac{1}{n^4}\Big) \cr
 \ar=\ar
\frac{m}{2} v_{n+1}\Delta^2 v_n + \frac{m}{m-1}(\Delta v_n+\Delta v_{n+1})\Delta v_n + o\Big(\frac{1}{n^4}\Big) \cr
 \ar=\ar
-\frac{m^2}{2}v_{n+1}^2\Delta v_n + \frac{m}{m-1}(\Delta v_n+\Delta v_{n+1})\Delta v_n + o\Big(\frac{1}{n^4}\Big) \cr
 \ar=\ar
\frac{4(m+1)}{m(m-1)}\frac{1}{n^4} + o\Big(\frac{1}{n^4}\Big).
 \eeqlb
By applying the Stolz--Ces\`aro theorem we deduce that
 \beqaa
\lim_{n\to \infty} n^3\Big(\Delta v_n + \frac{m}{2}v_nv_{n+1}\Big)
 =
\lim_{n\to \infty} \frac{\Delta\big(\Delta v_n + \frac{m}{2}v_nv_{n+1}\big)} {\frac{1}{(n+1)^3}-\frac{1}{n^3}} =
-\frac{4(m+1)}{3m(m-1)}.
 \eeqaa
This proves the desired result. \eproof

\blemma\label{th-v_n=(asym)} In the critical case, as $n\to \infty$ we have
 \beqa\label{v_n=(asym)}
v_n=\frac{2}{mn} - \frac{4(m+1)}{3m(m-1)}\frac{\log n}{n^2} + o\Big(\frac{\log n}{n^2}\Big).
 \eeqa
\elemma

\bproof By \eqref{(nv_n)(n^2(v_n-v_{n+1}))}, \eqref{eq:1/v_{n+1}-1/v_n-m/2=...} and \eqref{eq:v_n-v_{n+1}}, we have
 \beqnn
\Big(\frac{1}{v_{n+1}}-\frac{1}{v_n}\Big) - \frac{m}{2}
 =
\frac{m(m+1)}{3(m-1)n} + o(n^{-1}).
 \eeqnn
Then by the Stolz--Ces\`aro theorem,
 \beqnn
\lim_{n\to \infty} \frac{\big(\frac{1}{v_n} - \frac{1}{v_0}\big) - \frac{m}{2}n}{\log n}
 \ar=\ar
\lim_{n\to \infty} \frac{\big(\frac{1}{v_{n+1}} -\frac{1}{v_n}\big)- \frac{m}{2}}{\log\big(1+\frac{1}{n}\big)}  \cr
 \ar=\ar
\lim_{n\to \infty} \frac{\big(\frac{1}{v_{n+1}}-\frac{1}{v_n}\big)-\frac{m}{2}}{\frac{1}{n}}
 =
 \frac{m(m+1)}{3(m-1)}.
 \eeqnn
It follows that
 \beqlb\label{eq:1/v_n=...}
\frac{1}{v_n}
 =
\frac{mn}{2} + \frac{m(m+1)}{3(m-1)}\log n + o(\log n),
 \eeqlb
and hence
 \beqaa
v_n=
\frac{\frac{2}{m n}}{1 + \frac{2(m+1)}{3(m-1)}\frac{\log n}{n} + o\big(\frac{\log n}{n}\big)}.
 \eeqaa
By using Taylor's expansion for the function $1/(1+x)$, we get \eqref{v_n=(asym)}.  \eproof

\bremark\label{re-iterativescheme} We can give further refined asymptotic expansions for the sequence $\{v_n\}$ using an iterative argument. The approach goes as follows. Once a prior expansion for $v_n$ is obtained, we can substitute it into the r.h.s.\ of \eqref{eq:1/v_{n+1}-1/v_n-m/2=...} to get an expression for the difference $1/v_{n+1} - 1/v_n$. Then, by the Stolz--Ces\`aro theorem, we get an asymptotic expression for $1/v_n$, which yields a refined expansion for $v_n$ as in the proof of Lemma~\ref{th-v_n=(asym)}. This procedure can be repeated for further refinements. \eremark

\medskip\noindent\textit{Proof of Theorem~\ref{th-mixgeoasy0}-(3).~} In the critical case, the limits defined in \eqref{r_*,p_*=def} satisfy $r_*= 1-m^{-1}< 1$ and $p_*= 1$. Then \eqref{m^{-n}E(Y_n)=..} implies $F_\infty= \lim_{n\to \infty} m^{-n}\mbf{E}(Y_n)= 0$. The desired expansion for $r_n$ follows from \eqref{v_n=(asym)}. Since $v_n= O(n^{-1})$ and $v_n-v_{n+1}= O(n^{-2})$, from \eqref{eq:v_n-v_{n+1}} and \eqref{v_n=(asym)} it follows that
 \beqlb\label{v_n-v_{n+1}=..}
v_n-v_{n+1}\ar=\ar \frac{m}{2}v_nv_{n+1} + \frac{4(m+1)}{3m(m-1)}\frac{1}{n^3} + o\Big(\frac{1}{n^3}\Big) \cr
 \ar=\ar
\frac{m}{2}v_n^2 + \frac{m}{2}v_n(v_n-v_{n+1}) + \frac{4(m+1)}{3m(m-1)}\frac{1}{n^3} + o\Big(\frac{1}{n^3}\Big) \cr
 \ar=\ar
\frac{2}{mn^2} - \frac{8(m+1)}{3m(m-1)}\frac{\log n}{n^3} + o\Big(\frac{\log n}{n^3}\Big).
 \eeqlb
In view of \eqref{1-p_n=..}, we have
 \beqa\label{1-p_n=....}
1-p_n= \frac{r_n-r_{n+1}}{(m-1)r_{n+1}}
 =
\frac{m}{(m-1)^2} \frac{v_n-v_{n+1}}{1+v_{n+1}/(1-m^{-1})}.
 \eeqa
Then we can use \eqref{v_n=(asym)}, \eqref{v_n-v_{n+1}=..} and Taylor's expansion for the function $1/(1+x)$ to see
 \beqaa
1-p_n= \frac{m}{(m-1)^2}\Big[\frac{2}{mn^2} - \frac{8(m+1)}{3m(m-1)}\frac{\log n}{n^3} + o\Big(\frac{\log n}{n^3}\Big)\Big].
 \eeqaa
This proves the desired expansion of $p_n$. \qed

\medskip\noindent\textit{Proof of Corollary~\ref{th-sustprob0}.~} By Theorem \ref{th-mixgeo-mar0}, the random variable $Y_n$ has distribution $\mu_n= G(r_n,p_n)$. It follows that
 \beqnn
\mbf{P}(Y_n\ge 1)= {\blue 1 - p_n}.
 \eeqnn
{\blue The expression in the supercritical case follows immediately by Theorem~\ref{th-mixgeoasy0}}\qed

\medskip\noindent\textit{Proof of Corollary~\ref{th-condlimitlaw0}.~} (1)~Since $Y_n$ is distributed according to the geometric-type law $\mu_n= G(r_n,p_n)$, for any integer $M\ge 1$,
 \beqlb\label{P(liminfY_nleM)}
\mbf{P}\Big(\liminf_{n\to\infty} Y_n\le M\Big)
 \ar=\ar
\mbf{P}\Big(\bigcap_{n=1}^\infty \bigcup_{k=n}^\infty \{Y_n\le M\}\Big)
 \le
\lim_{n\to \infty} \sum_{k=n}^\infty \mbf{P}(Y_k\le M) \cr
 \ar=\ar
\lim_{n\to \infty}\sum_{k=n}^\infty \Big[p_k + (1-p_k)r_k\sum_{i=1}^M (1-r_k)^{i-1}\Big] \cr
 \ar=\ar
\lim_{n\to \infty}\sum_{k=n}^\infty \big\{p_k + (1-p_k)[1-(1-r_k)^M]\big\} \cr
 \ar\le\ar
\lim_{n\to \infty}\sum_{k=n}^\infty \big\{p_k + [1-(1-r_k)^M]\big\}.
 \eeqlb
In the supercritical case, by Theorem~\ref{th-mixgeoasy0}~(1) we have, as $k\to \infty$,
 \beqlb\label{r_kandp_k}
r_k= O(m^{-k}),
 \quad
p_k= O({\blue km^{-k}})
 \eeqlb
and
 \beqnn
(1-r_k)^M= 1 - Mr_k + o(r_k)= 1 + o(m^{-k}).
 \eeqnn
Then the series on the r.h.s.\ of \eqref{P(liminfY_nleM)} is convergent, implying $\mbf{P}(\liminf_{n\to \infty} Y_n\le M)= 0$, and so $Y_n\to \infty$ almost surely. Clearly, the Laplace transform of $m^{-n}Y_n$ is given by
 \beqnn
\mbf{E}(\e^{-s m^{-n}Y_n})
 =
p_n + \frac{(1-p_n)r_n\e^{-sm^{-n}}}{1 - \e^{-sm^{-n}} + r_n\e^{-sm^{-n}}} \cr
 =
p_n +\frac{ 1-p_n }{(\e^{sm^{-n}}-1)/r_n+1}.
 \eeqnn
By \eqref{r_kandp_k} we have
 \beqnn
\lim_{n\to \infty} \mbf{E}(\e^{-s m^{-n}Y_n})
 =
\frac{1}{F_{\infty}s+1}
 =
\frac{F_{\infty}^{-1}}{s+F_{\infty}^{-1}},
 \eeqnn
which implies the desired weak convergence result.

(2)~In the subcritical case, by Theorem~\ref{th-mixgeoasy0}~(2) we have $(1-p_n)= O(\gamma_*^{-n})$ as $n\to \infty$. It follows that
 \beqnn
\mbf{P}\Big(\limsup_{n\to \infty} Y_n\ge 1\Big)
 \ar=\ar
\mbf{P}\Big(\bigcap_{n=1}^\infty \bigcup_{k=n}^\infty \{Y_n\ge 1\}\Big)
 \leq
\lim_{n\to \infty} \sum_{k=n}^\infty \mbf{P}(Y_k\ge 1) \cr
 \ar=\ar
\lim_{n\to \infty} \sum_{k=n}^\infty (1-p_k)= 0,
 \eeqnn
which implies $Y_n\to 0$ almost surely. Moreover, since $\lim_{n\to \infty} r_n= r_*$, we have
 \beqnn
\lim_{n\to \infty}\mbf{E}(\e^{-sY_n}|Y_n\ge 1)
 =
\lim_{n\to \infty} \frac{r_n\e^{-s}}{1-(1-r_n)\e^{-s}}
 =
\frac{r_*\e^{-s}}{1-(1-r_*)\e^{-s}}.
 \eeqnn
This yields the desired weak convergence of the conditional distribution.

(3)~In the critical case, by Theorem~\ref{th-mixgeoasy0}~(3) we have $(1-p_k)= O(k^{-2})$ and $r_n\to 1-m^{-1}$ as $n\to \infty$. The desired results follow similarly as in the proof of part (2). \qed

\medskip\noindent\textit{Proof of Corollary~\ref{th-firmom0}.~} Since $Y_n$ has the geometric type distribution $\mu_n= G(r_n,p_n)$, we have  $\mbf{E}(Y_n)= (1-p_n)r_n^{-1}$. In the supercritical case, by Theorem~\ref{th-mixgeoasy0}~(1) and Taylor's expansion for the function $1/(1-x)$ we have
 \beqnn
\mbf{E}(Y_n)\ar=\ar \Big[1 - \frac{n}{F_\infty m^{n-1}} - \Big(\frac{p_0}{r_0} - mQ\Big)\frac{1}{F_\infty m^n} + \frac{n^2}{F_\infty^2m^{2(n-1)}} + o\Big(\frac{n^2}{m^{2n}}\Big)\Big] \cr
 \ar\ar
\cdot\,\Big[\frac{1}{F_\infty m^n} - \frac{n}{F_\infty^2m^{2n-1}} + o\Big(\frac{n}{m^{2n}}\Big)\Big]^{-1} \cr
 \ar=\ar
F_\infty m^n\Big[1 + {\blue o\Big(\frac{n}{m^{n}}\Big)}\Big].
 \eeqnn
In the subcritical case, by Theorem~\ref{th-mixgeoasy0}~(2) and Taylor's expansion for the function $1/(1+x)$,
 \beqnn
\mbf{E}(Y_n)\ar=\ar \Big[\frac{Kr_*\gamma_*^n}{m-1} + {\blue \Big(1 +\frac{mr_*} {1-\gamma_*}\Big) }\frac{{\blue K^2}r_*^2\gamma_*^{2n}}{(m-1)(1-\gamma_*)} + o(\gamma_*^{2n})\Big] \cr
 \ar\ar
\cdot\,\Big[r_* + \frac{Kr_*^2\gamma_*^n}{1-\gamma_*} + \Big(1+\frac{mr_*}{1-\gamma_*^2}\Big) \frac{{\blue K^2}r_*^3\gamma_*^{2n}}{(1-\gamma_*)^2} + o(\gamma_*^{2n})\Big]^{-1} \cr
 \ar=\ar
\frac{K\gamma_*^n}{m-1} + {\blue \frac{mK^2r_*\gamma_*^{2n}}{(m-1)(1-\gamma_*)^2} } + o(\gamma_*^{2n}).
 \eeqnn
In the critical case, by Theorem~\ref{th-mixgeoasy0}~(3) and Taylor's expansion,
 \beqnn
\mbf{E}(Y_n)\ar=\ar \Big[\frac{2}{(m-1)^2n^2} + \frac{8(m+1)\log n}{3(m-1)^3n^3} + o\Big(\frac{\log n}{n^3}\Big)\Big] \cr
 \ar\ar
\cdot\,\Big[\frac{m-1}{m} + \frac{2}{mn} - \frac{4(m+1)\log n}{3m(m-1)n^2} + o\Big(\frac{\log n}{n^2}\Big)\Big]^{-1} \cr
 \ar=\ar
\frac{m}{m-1}\Big[\frac{2}{(m-1)^2n^2} + \frac{8(m+1)\log n}{3(m-1)^3n^3} + o\Big(\frac{\log n}{n^3}\Big)\Big].
 \eeqnn
Then we have the desired expansions in the three cases. \qed

\medskip\noindent\textit{Proof of Corollary~\ref{th-E(.^{Y_n})-1=..0}.~} By Theorem~\ref{th-mixgeo-mar0}, the random variable $Y_n$ has distribution $\mu_n= G(r_n,p_n)$, so its generating function is given by \eqref{E(s^{Y_n})=p_n+..}.

(1)~In the supercritical case, by Taylor's expansion for the function $1/(1+x)$, for any $|s|< 1$ we have
 \beqnn
\mbf{E}(s^{Y_n})\ar=\ar p_n + \frac{(1-p_n)r_ns}{1-s} \frac{1}{1 + r_ns/(1-s)} \cr
 \ar=\ar
p_n + \frac{(1-p_n)r_ns}{1-s} \Big[1 - \frac{r_ns}{1-s} + o(r_n)\Big] \cr
 \ar=\ar
p_n + \frac{r_ns}{1-s} + o(r_np_n).
 \eeqnn
Then we use Theorem~\ref{th-mixgeoasy0}-(1) to see that
 \beqnn
\mbf{E}(s^{Y_n})
 \ar=\ar
\frac{n}{F_\infty m^{n-1}} + \Big(\frac{p_0}{r_0} - mQ\Big)\frac{1}{F_\infty m^n} - \frac{n^2}{F_\infty^2m^{2(n-1)}} + o\Big(\frac{n^2}{m^{2n}}\Big) \cr
 \ar\ar
+\, \frac{s}{1-s}\Big[\frac{1}{F_\infty m^n} - \frac{n}{F_\infty^2m^{2n-1}} + o\Big(\frac{n}{m^{2n}}\Big)\Big].
 \eeqnn
By ignoring the higher order infinitesimal terms, we get the desired expansion.

(2)~In the subcritical case, by \eqref{E(s^{Y_n})=p_n+..} we have, for $0< s< (1-r_n)^{-1}$,
 \beqnn
\mbf{E}(s^{Y_n})-1= (1-p_n) \Big[\frac{r_ns}{1-(1-r_n)s} - 1\Big]
 =
\frac{(1-p_n)(s-1)}{1-(1-r_n)s},
 \eeqnn
where
 \beqnn
1-(1-r_n)s= [1-(1-r_*)s]\Big[1 + \frac{(r_n-r_*)s}{[1-(1-r_*)s]}\Big].
 \eeqnn
Then, by Taylor's expansion and Theorem~\ref{th-mixgeoasy0}-(2),
 \beqlb\label{E(s^{Y_n})-1=1+..(0)}
\mbf{E}(s^{Y_n})-1\ar=\ar \frac{(1-p_n)(s-1)}{1-(1-r_*)s} \Big[1 + \frac{(r_n-r_*)s}{[1-(1-r_*)s]}\Big]^{-1} \cr
 \ar=\ar
\frac{s-1}{1-(1-r_*)s} \Big\{\frac{Kr_*\gamma_*^n}{m-1} + {\blue \Big(1 + \frac{mr_*} {1-\gamma_*}\Big) }\frac{{\blue K^2}r_*^2\gamma_*^{2n}}{(m-1)(1-\gamma_*)} \cr
 \ar\ar
+\, o(\gamma_*^{2n})\Big\} \cdot\Big\{1 - \frac{s}{1-(1-r_*)s} \Big[\frac{Kr_*^2\gamma_*^n}{1-\gamma_*} + \Big(1+\frac{mr_*}{1-\gamma_*^2}\Big) \frac{{\blue K^2}r_*^3\gamma_*^{2n}}{(1-\gamma_*)^2}\Big] \cr
 \ar\ar
+\, \frac{s^2}{[1-(1-r_*)s]^2}\frac{{\blue K^4}r_*^4\gamma_*^{2n}}{(1-\gamma_*)^2} + o(\gamma_*^{2n})\Big\}.
 \eeqlb
This leads to the desired expansion of the generating function.

(3)~In the critical case, we have $r_*= 1-m^{-1}$. As in \eqref{E(s^{Y_n})-1=1+..(0)} one can use Theorem~\ref{th-mixgeoasy0}-(3) and Taylor's expansion to see that, for $0< s< m$,
 \beqnn
\mbf{E}(s^{Y_n})-1\ar=\ar \frac{(1-p_n)(s-1)}{1-m^{-1}s} \Big[1 + \frac{(r_n-r_*)s}{1-m^{-1}s}\Big]^{-1} \cr
 \ar=\ar
\frac{s-1}{1-m^{-1}s} \Big[\frac{2}{(m-1)^2n^2} + \frac{8(m+1)\log n}{3(m-1)^3 n^3} + o\Big(\frac{\log n}{n^3}\Big)\Big] \cr
 \ar\ar
\cdot\,\Big\{1 - \frac{s}{1-m^{-1}s} \Big[\frac{2}{mn} - \frac{4(m+1)\log n}{3m(m-1)n^2}\Big] + o\Big(\frac{1}{n^2}\Big) \cr
 \ar\ar
+\, \frac{s^2}{[1-m^{-1}s]^2}\frac{4}{m^2n^2}\Big\}.
 \eeqnn
This gives the expansion \eqref{E(s^{Y_n})=1+..}. By Theorem~\ref{th-mixgeoasy0}-(3) one can also see that
 \beqaa
1-(1-r_n)m=\frac{2}{n}\Big[1 - \frac{2(m+1)\log n}{3(m-1)n} + o\Big(\frac{\log n}{n}\Big)\Big].
 \eeqaa
As in \eqref{E(s^{Y_n})-1=1+..(0)}, by Taylor's expansion we have
 \beqnn
\mbf{E}(m^{Y_n})-1\ar=\ar \frac{(1-p_n)(m-1)n}{2}\Big[1 - \frac{2(m+1)\log n}{3(m-1)n} + o\Big(\frac{\log n}{n}\Big)\Big]^{-1} \cr
 \ar=\ar
(m-1)n\Big[\frac{1}{(m-1)^2n^2} + \frac{4(m+1)\log n}{3(m-1)^3 n^3} + o\Big(\frac{\log n}{n^3}\Big) \cr
 \ar\ar
+\,\frac{2(m+1)\log n}{3(m-1)^3n^3}\Big].
 \eeqnn
This proves the desired expansion \eqref{E(m^{Y_n})=1+..}. \qed

\bigskip

\textbf{Acknowledgements.~} We would like to thank Professors Xinxing Chen and Yueyun Hu for enlightening discussions. This research is supported by the National Key R{\&}D Program of China (No.~2020YFA0712901).

\textbf{Data Availability Statement.~} Data sharing not applicable to this article as no datasets were generated or analyzed during the current study.

\textbf{Conflict of Interest Statement.~} The authors declare that they have no conflict of interest.

\end{document}